\documentstyle[12pt,amssymb,amstex]{amsart}
\def\cb{{\cal B}}

\def\cf{{\cal F}}
\def\cg{{\cal G}}
\def\ch{{\cal H}}

\def\car{{\cal R}}
\def\cs{{\cal S}}

\def\ga{{\frak A}}

\def\bp{{\mathbb P}}

\def\bz{{\mathbb Z}}

\def\a{\alpha}
\def\b{\beta}
  
\def\d{\delta}

\def\k{\kappa}
\def\m{\mu}

\def\r{\rho}
\def\s{\sigma} 

\def\f{\varphi} 
\def\th{\theta}  
\def\om{\omega} \def\Om{\Omega}

\newtheorem{thm}{Theorem}%[section]
\newtheorem{lem}[thm]{Lemma}
\newtheorem{cor}[thm]{Corollary}
\newtheorem{prop}[thm]{Proposition}

\newtheorem{rem}[thm]{Remark}

\def\di{\mathop{\rm d}}

\def\idd{{\bf 1}\!\!{\rm I}}

\begin{document}

\title[$q$--relations, uniquely mixing]
{unique mixing of the shift on the $C^*$--algebras generated by the $q$--canonical commutation relations}
\author[Dykema]{Kenneth Dykema$^*$}\thanks{${}^*$Research supported in part by NSF grant DMS--0600814}
\address{Kenneth Dykema\\
Department of Mathematics\\
Texas A\&M University \\
College Station TX 77843--3368, USA}
\email{{\tt kdykema@@math.tamu.edu}}
\author[Fidaleo]{Francesco Fidaleo${}^\dag$}\thanks{${}^\dag$Permanent address:
Dipartimento di Matematica,
II Universit\`{a} di Roma ``Tor Vergata'',
Via della Ricerca Scientifica, 00133 Roma, Italy.
e--mail: \tt fidaleo@@mat.uniroma2.it}
\address{Francesco Fidaleo\\
Department of Mathematics\\
Texas A\&M University \\
College Station TX 77843--3368, USA}

\begin{abstract}
The shift on the $C^*$--algebras generated by the Fock representation of the $q$--commutation
relations has the strong ergodic property of unique mixing, when $|q|<1$.
\end{abstract}

\maketitle

\section{introduction}
The $q$--commutation relations have been studied in the physics literature, see e.g. \cite{FB}.
These are the relations
\[
a_ia_j^+-qa_j^+a_i=\delta_{ij}1,\qquad i,j\in\mathbb{Z}
\]
where $-1\le q\le 1$.
This gives an interpolation between the canonical commutation relations (Bosons) when $q=1$ and the canonical anticommutation
relations (Fermions) when $q=-1$, while when $q=0$ we have freeness (cf. \cite{VDN}).
In~\cite{BS1}, (see also \cite{Gr} and \cite{Fiv})
a Fock representation of these relations was found,
giving annihilators $a_i$ and their adjoints, the creators $a_i^+$, acting on a Hilbert space with
a vacuum vector $\Omega$.
The $C^*$--algebras and von Neumann algebras generated by sets of these operators or by their real parts $a_i+a_i^+$ have been
much studied. The reader is referred to \cite{JSW, DN, Nou, Shl, Sn, Ri}
for  results, applications and further details. 

In the present note, we show that the shift $\alpha_q$ on the C$^*$--algebras generated by these $a_i$, or by their self--adjoint parts,
has the strong ergodic property of unique mixing, which was introduced in the companion paper~\cite{F}, whenever $|q|<1$.
The case of free shifts was treated in~\cite{F}, but
when $q\ne0$, this shift $\alpha_q$ cannot a free shift,
so these results provide new noncommutative examples of unique mixing.
We also observe that the examples provided by the $a_i$ and by the $a_i+a_i^+$
are nonconjugate.
It was shown in \cite{F} that there is no classical
counterpart to this situation.

\section{terminology and basic notations}

For a (discrete) $C^*$--dynamical system we mean a triplet $\big(\ga,\a,\om\big)$ 
consisting of a unital $C^*$-algebra $\ga$, an automorphism $\a$ of $\ga$, 
and a state $\om\in\cs(\ga)$ invariant under the action of $\a$. The pair $(\ga,\a)$ consisting of 
$C^*$-algebra and an automorphism as before, is called a $C^*$--dynamical system as well. A classical
$C^*$--dynamical system is simply a dynamical system such that $\ga\sim C(X)$, $C(X)$ being the Abelian $C^*$-algebra of all the continuous functions on the compact Hausdorff space $X$. In this situation, $\a(f)=f\circ T$ for some homeomorphism $T:X\mapsto X$. 

Consider for $j=1,2$, the $C^*$--dynamical systems $(\ga_j,\a_j,\om_j)$ together with the canonically associated $W^*$--dynamical systems $(M_j,\hat\a_j,\hat\om_j)$. Here, $M_j:=\pi_{\om_j}(\ga_j)''$ and 
$\hat\a_j$, $\hat\om_j $ are the canonical extensions of $\a_j$ and $\om_j$ to $M_j$, respectively.
The $C^*$--dynamical system $(\ga_j,\a_j,\om_j)$ are said to be {\it conjugate} if there exists an automorphism $\b:M_1\mapsto M_2$ intertwining the dynamics ($\b\hat\a_1=\hat\a_2\b$), and the states ($\hat\om_2\circ\b=\hat\om_1$). Suppose that $\ga_j\sim C(X_j)$, $X_j$ being compact spaces.  Then there exist probability measures $\m_j$ on $X_j$, and measure-preserving 
homeomorphisms $T_j$ of the compact spaces $X_j$ such that
$$
\a_j(f)=f\circ T_j\,,\qquad \om_j(f)=\int_{X_j}f\di\m_j\,.
$$

Thanks to a result by J. von Neumann (cf. \cite{B}, p.\ 69), our definition is equivalent to the following one, provided that the $X_j$ are compact metric spaces. There exist 
$\m_j$--measurable sets $A_j\in X_j$ of full measure such that 
$T_j(A_j)=A_j$, and a one--to--one measure--preserving map $S:A_1\mapsto A_2$ such that
$T_2=S\circ T_1\circ S^{-1}$. The reader is referred to \cite{Je} for further details relative to the classical case. 

To recall the definition from~\cite{F},
a $C^*$--dynamical system $(\ga,\a)$ is said to be
{\it uniquely mixing} if
\begin{equation*}
%\label{mppp1}
\lim_{n\to+\infty}\f(\a^{n}(x))=\f(\idd)\om(x)\,,\quad x\in\ga\,,\f\in\ga^*
\end{equation*}
for some $\om\in\cs(\ga)$.

It can readily seen that $\om$ is invariant under $\a$. In addition, it is unique among the invariant states for $\a$. 

\bigskip

Let $\ch:=\ell^2(\bz)$, with $e_i\in\ch$ the function taking value $1$ at $i$ and zero elsewhere.
The $q$--deformed Fock space $\cf_q$ is the completion of the algebraic linear span of the vacuum vector $\Om$, together with vectors
$$
f_1\otimes\cdots\otimes f_n\,,\quad f_j\in\ch\,,j=1,\dots,n\,,n=1,2,\dots
$$
with respect to the inner product
$$
\langle f_1\otimes\cdots\otimes f_n\,,g_1\otimes\cdots\otimes g_m\rangle_q
:=\d_{n,m}\sum_{\pi\in\bp_n}q^{i(\pi)}\langle f_1\,,g_{\pi(1)}\rangle\cdots\langle f_n\,,g_{\pi(n)}\rangle\,,
$$
$\bp_n$ being the symmetric group of $n$ elements, and $i(\pi)$ the number of inversions of $\pi\in\bp_n$. We have $\langle f\,,g\rangle_q=\langle P_qf\,,g\rangle_0$, where $P_q$ is determined by
\begin{equation}
\label{pii}
P_q\Om=\Om\,,\quad P_q f_1\otimes\cdots \otimes f_n=
q^{i(\pi)}f_{\pi(1)}\otimes\cdots \otimes f_{\pi(n)}\,.
\end{equation}
The creator $a_i^+$ acts on $\cf_q$ by 
$$
a_i^+\Omega=e_i\,,\quad
a_i^+(f_1\otimes\cdots\otimes f_n)=e_i\otimes f_1\otimes\cdots\otimes f_n\,,
$$
and its adjoint is the annihilator $a_i$ given by 
\begin{align*}
a_i\Omega&=0\,,\\
a_i(f_1\otimes\cdots\otimes f_n)&=\sum_{k=1}^nq^{k-1}\langle f_k,e_i\rangle
f_1\otimes\cdots f_{k-1}\otimes f_{k+1}\otimes\cdots\otimes f_n\,.
\end{align*}

Denote by $\car_q$ the $C^{*}$--algebra generated by $\{a_i\mid i\in\mathbb{Z}\}$,
and by $\cg_q$
the $C^{*}$--algebra generated by $\{a_i+a^{+}_i\mid i\in\mathbb{Z}\}$.
The right shift $\a=\a_q$ acting on $\car_q$ is uniquely determined by
$$
\a(a_i):=a_{i+1}\,,\qquad i\in\bz
$$
on the generators.
The Fock vacuum expectation $\om:=\langle\,{\bf\cdot}\,\Om\,,\Om\rangle$ is invariant for the shift $\a$.
The restriction of the vacuum expectation to $\cg_q$ is a faithful trace. For further details, we refer
to \cite{BKS, BS} and the literature cited therein. 

\section{unique mixing of the q--shifts}

In the present section we prove the announced result on the unique mixing of the shift on the $q$--canonical commutation relations. We start with the following
\begin{lem}
\label{bo}
Let $\{\xi_j\}_{j=1}^{n}\subset\ch^{\otimes n}$, and $\{f_j\}_{j=1}^{n}\subset\ch$ be an orthonormal set. Then
$$
\bigg\|\sum_{j=1}^{n}a^{+}(f_j)\xi_j\bigg\|\leq\sqrt{\frac{n}{1-|q|}}\max_{1\leq j\leq n}\|\xi_j\|\,.
$$
\end{lem}
\begin{pf}
Denote as in \cite{BS}, $P^{(n)}_q:=P_q\lceil_{\ch^{\otimes n}}$, where $P_q$ is given in \eqref{pii}.
By taking into account Section 3 of \cite{BS} (see also \cite{BKS}, Section 1), we get
\begin{align*}
\bigg\langle\sum_{j=1}^{n}&a^{+}(f_j)\xi_j\,,\sum_{j=1}^{n}a^{+}(f_j)\xi_j\bigg\rangle_q
=\bigg\langle P_q^{(k+1)}\sum_{j=1}^{n}f_j\otimes\xi_j\,,\sum_{j=1}^{n}f_j\otimes\xi_j\bigg\rangle_0\\
\leq&\frac{1}{1-|q|}\bigg\langle\idd\otimes P_q^{(k+1)}\sum_{j=1}^{n}f_j\otimes\xi_j\,,\sum_{j=1}^{n}f_j\otimes\xi_j\bigg\rangle_0
  \displaybreak[2] \\
\leq&\frac{1}{1-|q|}\sum_{i,j=1}^{n}\big\langle f_i\otimes P_q^{(k)}\xi_i\,,f_j\otimes\xi_j\big\rangle_0 \displaybreak[3] \\
=&\frac{1}{1-|q|}\sum_{i,j=1}^{n}\langle f_i\,,f_j\rangle\langle P_q^{(k)}\xi_i\,,\xi_j\rangle_0 \displaybreak[3] \\
=&\frac{1}{1-|q|}\sum_{i,j=1}^{n}\langle f_i\,,f_j\rangle\langle\xi_i\,,\xi_j\rangle_q \displaybreak[1] \\
=&\frac{1}{1-|q|}\sum_{i=1}^{n}\langle\xi_i\,,\xi_i\rangle_q\\
\leq&\frac{n}{1-|q|}\max_{1\leq j\leq n}\|\xi_j\|^2\,.
\end{align*}
\end{pf}
\begin{prop}
\label{bo1}
Let $0\leq k_1<k_2<\cdots<k_n<\cdots$ be a sequence of increasing natural numbers, and
$e_{\s_1},\dots,e_{\s_i},e_{\r_1},\dots,e_{\r_j}$ elements of the canonical basis of $\ell^2(\bz)$. We have
$$
\bigg\|\sum_{l=1}^{n}\a^{k_l}\big(a^{+}(e_{\s_1})\cdots a^{+}(e_{\s_i})a(e_{\r_1})\cdots a(e_{\r_j})\big)
\bigg\|\leq\sqrt{\frac{n}{(1-|q|)^{i+j}}}
$$
if at least either $i$ or $j$ is nonnull.
\end{prop}
\begin{pf}
Suppose first $i>0$. It is enough consider unit vectors $\xi\in\ch^{\otimes m}$, $m=j,j+1,\dots$\,. Put
$$
\xi_l:=a^{+}(e_{\s_2+k_l})\cdots a^{+}(e_{\s_i+k_l})a(e_{\r_1+k_l})\cdots a(e_{\r_j+k_l})\xi\,.
$$
Notice that, by Theorem 3.1 of \cite{BKS}, $\|\xi_l\|\leq1/\sqrt{(1-|q|)^{i+j-1}}$. In addition,
$\langle e_{\s_1+k_l}\,,e_{\s_1+k_{\hat l}}\rangle=\d_{l,\hat l}$.
By applying Lemma \ref{bo}, we get
\begin{align*}
&\bigg\|\sum_{l=1}^{n}\a^{k_l}\big(a^{+}(e_{\s_1})\cdots a^{+}(e_{\s_i})a(e_{\r_1})
\cdots a(e_{\r_j})\big)\xi\bigg\|^2\\
=&\bigg\langle\sum_{l=1}^{n}a^{+}(e_{\s_1+k_l})\xi_l\,,\sum_{l=1}^{n}a^{+}(e_{\s_1+k_l})\xi_l
\bigg\rangle\leq\frac{n}{(1-|q|)^{i+j}}\,.
\end{align*}
If $i=0$ and then $j\neq0$ (i.e. we have only annihilators), the assertion follows by the first part as
$$
\sum_{l=1}^{n}\a^{k_l}\big(a(e_{\r_1})\cdots a^{+}(e_{\r_j})\big)
=\bigg(\sum_{l=1}^{n}a^+(e_{\r_j+\k_l})\cdots a^{+}(e_{\r_1+k_l})\bigg)^*\,.
$$
\end{pf}

The following theorem is the annouced result on the strong ergodic property enjoyned by the $q$--shift.
\begin{thm}
\label{bo2}
The dynamical system $(\car_q,\a)$ is uniquely mixing, with the vacuum expectation $\om$ as the unique invariant state.
\end{thm}
\begin{pf}
Let $X\in\car_q$ have vanishing vacuum expectation. It is norm limit of elements as those treated in Proposition \ref{bo1}. By Proposition 2.3 of~\cite{F}, and a standard approximation argument, it is enough to prove that
$$
\frac{1}{n}\bigg\|\sum_{l=1}^{n}\a^{k_l}(X)\bigg\|\longrightarrow0
$$
for each $X$ given by
$$
X= a^{+}(e_{\s_1})\cdots a^{+}(e_{\s_i})a(e_{\r_1})\cdots a(e_{\r_j})
$$
for which either $i$ or $j$ is nonnull, and for each sequence $0\leq k_1<k_2<\cdots<k_n<\cdots$ of increasing natural numbers. The result directly follows by Proposition \ref{bo1}.
\end{pf}

The restriction of any uniquely mixing automorphism to any invariant $C^*$--subalgebra is also uniquely mixing. For example,
\begin{cor} \label{umcor}
The dynamical system $(\cg_q,\a\lceil_{\cg_q})$ is uniquely mixing, with the vacuum expectation
$\om\lceil_{\cg_q}$ as the unique invariant state.
\end{cor}
\begin{rem}
The dynamical systems $(\car_q,\a^{-1})$, $(\cg_q,\big(\a\lceil_{\cg_q}\big)^{-1})$ are uniquely mixing as well, with the vacuum expectation as the unique invariant state.
\end{rem}
This can be shown by taking into account that
$$
\th\a=\a^{-1}\th\,,\qquad\om\circ\th=\om\,,
$$
where $\th(a(e_k)):=a(e_{-k})$, $k\in\bz$, is the "time reversal", $a(e_k)$ being the $k$--annihilator.

It is known from~\cite{BKS} that $\cg_q''$ is a II$_1$--factor.
Moreover, it is well known and easily seen that $\car_q''$ is all of $\cb(\cf_q)$.
(A proof in the case of the Fock representation of finitely many $a_i$ is found in~\cite{DN}).
For convenience, we provide a proof of this fact below.
In any case, from this it follows that $(\car_q,\alpha)$ is not conjugate to $(\cg_q,\alpha\lceil_{\cg_q})$
\begin{prop}
For all $-1<q<1$, the von Neumann algebra
$\car_q''$ is all of $\cb(\cf_q)$.
\end{prop}
\begin{pf}
It will suffice to show that the rank--one projection $P_\Omega$ onto the span of the vacuum vector
belongs to $\car_q''$, because $\Omega$ is cyclic for the action of $\car_q$ on $\cf_q$.
In the case of $q=0$, let us write $v_i$ for $a_i$ acting on $\cf_0$.  Thus, $\{v_i^*\}_{i\in\mathbb{Z}}$
is a family of isometries with orthogonal ranges, and
the sum $\sum_{i\in\mathbb{Z}}v_i^*v_i$ converges in strong operator topology
to $I-P_\Omega$.
For the general case, we fix $-1<q<1$ and we
will refer to Propositions 3.2 and 3.4 and Remark 3.3 of~\cite{DN}.
These show that there is a unitary $U:\cf_q\to\cf_0$ that sends the $n$--particle space to the $n$--particle
space, and there is a positive operator $M$ on $\cf_q$ given by $M\Omega=\Omega$ and
\[
M(f_1\otimes\cdots\otimes f_n)=\sum_{k=1}^nq^{k-1}f_k\otimes f_1\otimes\cdots f_{k-1}\otimes f_{k+1}\otimes\cdots\otimes f_n\,.
\]
Proposition 3.2 of~\cite{DN} shows that $M$ has range equal to all of $\cf_q$.
Moreover, we have $Ua_iU^*=v_iR$,
where $R=UM^{1/2}U^*$.
Thus,
\[
U(\sum_{i\in\mathbb{Z}}a_i^+a_i)U^*=R(\sum_{i\in\mathbb{Z}}v_i^*v_i)R=R(I-P_\Omega)R,
\]
where the sums converge in strong operator topology.
From this, we see that the range projection of $\sum_{i\in\mathbb{Z}}a_i^+a_i$ is $I-P_\Omega$.
\end{pf}

Notice that $(\car_q,\a,\om)$ is not conjugate to $(\cg_q,\a\lceil_{\cg_q},\om\lceil_{\cg_q})$. Indeed, $\pi_{\om}(\car_q)''\equiv\car_q'' $ is not isomorphic to $\pi_{\om}(\cg_q)''\equiv\cg_q'' $ as we have shown.  Thus, we provide nontrivial examples of uniquely mixing $C^*$--dynamical systems for which  the unique invariant state is faithful (case of the $C^{*}$--algebra $\cg_q$ generated by the self--adjoint part for which the restriction of vacuum state is a faithful trace), or when it is not faithful (case of the $C^{*}$--algebra $\car_q$). However, for all cases the associated GNS covariant representation is faithful. It was shown in \cite{F} that there is no classical counterpart to this situation.

\end{document}